\documentclass[12pt,a4paper]{amsart}
\usepackage[latin1]{inputenc}
\usepackage{latexsym}
\usepackage{color,graphicx,shortvrb}
\usepackage{amsmath, amssymb}
\usepackage{amsfonts}
\usepackage[colorlinks, bookmarks=true]{hyperref}

\newtheorem{theorem}{Theorem}[section]
\newtheorem{lemma}[theorem]{Lemma}

\theoremstyle{definition}
\newtheorem{definition}[theorem]{Definition}
\newtheorem{remark}[theorem]{Remark}

\author{J. M. Almira}
\title{An Ultrametric lethargy result and its application to $p$-adic number theory}

\begin{document}

\begin{abstract} In this paper we show a lethargy result in the non-Archimedean context, for general ultrametric approximation schemes and, as a consequence, we  prove the existence of $p$-adic transcendental numbers whose best approximation errors by algebraic $p$-adic numbers of degree $\leq n$ decays slowly.
\end{abstract}

\maketitle

\markboth{J. M. Almira}{An Ultrametric lethargy result and its applications}

\section{Motivation}

Let $(X,\|\cdot\|)$ be an ultrametric Banach space over a non-Archimedean valued field $\mathbb{K}$, and let
$\{0\}=A_0\subset A_1\subset\ldots \subset A_n\subset\ldots \subset X$
be an infinite chain of subsets of $X$, where all inclusions are
strict.

\begin{definition} We say that $(X,\{A_n\})$ is an {\it ultrametric approximation scheme} (or that $(A_n)$ is an approximation scheme in $X$) if:
\begin{itemize}
\item[$(i)$] There exists a map $K:\mathbb{N}\to\mathbb{N}$ such that $K(n)\geq n$ and $A_n+A_n\subseteq A_{K(n)}$ for all $n\in\mathbb{N}$.

\item[$(ii)$] $\lambda A_n\subset A_n$ for all $n\in\mathbb{N}$ and all scalars $\lambda\in \mathbb{K}$.

\item[$(iii)$] $\bigcup_{n\in\mathbb{N}}A_n$ is a dense subset of $X$
\end{itemize}
The approximation scheme is called {\it non-trivial} if $X \neq \cup_n \overline{A_n}$.
\end{definition}

A particular example is a {\it linear approximation scheme}, arising when the
sets $A_n$ are linear subspaces of $X$.  Approximation schemes were introduced, for $X$ a quasi-Banach space over $\mathbb{R}$ or $\mathbb{C}$, by Butzer and Scherer in 1968 \cite{butzer_scherer}, and are a natural object to study in classical approximation theory. Their study leads to the theory of approximation spaces and, in particular, the problem of characterizing membership to several  functional approximation spaces in terms of smoothness. To be more precise, the so called Central Theorems in approximation theory show that, in general, there exists a strong relation between the rate of decay of the best approximation errors $E(f,A_n)$ and the smoothness properties of $f$ \cite{devore,LGM,Pie}.  On the other hand, one of the most remarkable early results in the constructive theory of functions
is Bernstein Lethargy Theorem, which states that if
$X_0\subsetneq X_1\subsetneq X_2\subsetneq \cdots \subsetneq X$ is an ascending chain of finite dimensional vector subspaces of a Banach space $X$, and $\{\varepsilon_n\}\searrow 0$
is a non-increasing sequence of positive real numbers that converges to zero, then there exists an element $x\in X$ such that the $n$-th error of best approximation by elements
of $X_n$ satisfies $E(x,X_n)=\varepsilon_n$ for all $n\in\mathbb{N}$.   This result was first obtained in 1938 by  Bernstein \cite{bernsteininverso}
for $X=C([0,1])$ and $X_n=\Pi_n$, the vector space of real polynomials of degree $\leq n$.
The case of arbitrary finite dimensional $X_n$ is treated, for instance, in
\cite[Section II.5.3]{singerlibro}.

There are several generalizations of Bernstein's result to arbitrary chains of
(possibly infinite dimensional) closed subspaces
$X_1 \subsetneq X_2 \subsetneq \ldots$ of a Banach space $X$.
For example, Tjuriemskih \cite{tjuriemskih} and Nikolskii \cite{nikolskii,nikolskii2}
(see also \cite[Section I.6.3]{singerlibro}) proved that a sufficient (resp.~necessary)
condition for the existence of $x \in X$ verifying $E(x,X_n) = \varepsilon_n$ is that
$X$ is a Hilbert space (resp.~$X$ is reflexive). These results were proved independently and
by other means by Almira and Luther \cite{almiraluther1, almiraluther2} and Almira and Del Toro \cite{almiradeltoro1,almiradeltoro2}. Also, Bernstein
Lethargy Theorem has been generalized to
chains of finite-dimensional subspaces in non-Banach spaces (such as $SF$-spaces)
by  Lewicki \cite{lewicki, lewicki1}. These two approaches were successfully combined
by Micherda \cite{micherda}. Finally, thanks to the work by Plesniak \cite{plesniak}, the lethargy theorem has become a very
useful tool for the theory of quasianalytic functions of several complex variables.

In 1964 Shapiro \cite{shapiro} used Baire Category Theorem to prove that,
for any sequence $X_1 \subsetneq X_2 \subsetneq \ldots \subsetneq X$ of
closed (not necessarily finite dimensional) subspaces of a Banach space $X$,
and any sequence $\{\varepsilon_{n}\}\searrow 0$, there exists an $x\in X$ such that $E(x,X_{n})\neq\mathbf{O}(\varepsilon_{n})$.  This result was strengthened by Tjuriemskih \cite{tjuriemskih1} who, under the very same conditions of Shapiro's Theorem, proved the existence of  $x\in X$ such that $E(x,X_{n})\geq \varepsilon_{n}$, $n=0,1,2,\cdots$. Moreover, Borodin \cite{borodin} gave an easy proof of this result and proved that, for arbitrary infinite dimensional Banach spaces $X$ and for sequences $\{\varepsilon_n\}\searrow 0$ satisfying $\varepsilon_n>\sum_{k=n+1}^\infty\varepsilon_k$, $n=0,1,2,\cdots$, there exists  $x\in X$ such that $E(x,X_{n})= \varepsilon_{n}$, $n=0,1,2,\cdots$.

In this paper we show a result analogous to Shapiro's Theorem, but in the non-Archimedean context and for general ultrametric approximation schemes (Theorem \ref{main}) and, as a consequence, we prove the existence of $p$-adic transcendental numbers whose best approximation errors by algebraic p-adic numbers of degree $\leq n$ decays slowly (Theorem \ref{transcendental}).

Throughout the paper we consider $(X,\|\cdot\|)$ an ultrametric Banach space over a non-Archimedean complete valued field $(\mathbb{K},|\cdot|)$ which is nontrivial.  By nontrivial, we mean that the absolute value $|\cdot|$ is not the trivial one, which is given by $|x|=1$ for $x\neq 0$ and $|0|=0$. If $\mathbb{K}$ is such a field and
\[
p_{\mathbb{K}}=\sup\{|\lambda|: \lambda\in\mathbb{K}, |\lambda|<1\}, \]
then either $p_{\mathbb{K}}<1$, in which case, the set $|\mathbb{K}^{\times}|:=\{|\lambda|:\lambda\in\mathbb{K}\setminus\{0\}\}$ satisfies  $|\mathbb{K}^{\times}|=\{p_{\mathbb{K}}^n:n\in\mathbb{Z}\}$ and there exists $\rho\in \mathbb{K}$ such that $p_{\mathbb{K}}=|\rho|$, or $p_{\mathbb{K}}=1$, in which case,  $|\mathbb{K}^{\times}|$ is a dense subset of $(0,\infty)$. If $p_{\mathbb{K}}<1$ we say that the absolute value of $\mathbb{K}$ is discrete. If $p_{\mathbb{K}}=1$, we say that $\mathbb{K}$ is densely valued.  The most important examples of non-Archimedean complete valued fields we deal with are the fields of $p$-adic numbers $\mathbb{Q}_p$ and $\mathbb{C}_p$. Concretely, the elements of $\mathbb{Q}_p$ are the expressions of the form
\begin{equation}\label{numero}
x=a_{m}p^{m}+a_{m+1}p^{m+1}+\cdots+a_0+a_1p+a_2p^2+ \cdots + a_np^n +\cdots =\sum_{n\geq m}a_np^n,
\end{equation}
where $m\in\mathbb{Z}$, $1\leq a_{m}\leq p-1$ and  $0\leq a_k\leq p-1$ for all $k>m$. Given $x$ as in \eqref{numero},  its $p$-adic absolute value is given by $|x|_p=p^{-m}$. Obviously, the set of rational numbers $\mathbb{Q}$ is a dense subset of $\mathbb{Q}_p$. Indeed, $\mathbb{Q}_p$ results from $\mathbb{Q}$ by topological completion, when we consider over $\mathbb{Q}$ the topology defined by the p-adic absolute value $|\cdot|_p$. The construction of $\mathbb{C}_p$ is a little bit more difficult. As a first step, it is shown that if $\mathbb{K}$ is an algebraic extension of $\mathbb{Q}_p$, there exists a unique absolute value over $\mathbb{K}$ extending the absolute value of $\mathbb{Q}_p$. This implies that there exists a unique absolute value over $\mathbb{Q}_p^a$, the field of algebraic numbers over $\mathbb{Q}_p$, which extends the p-adic absolute value $|\cdot|_p$. Then $\mathbb{C}_p$ results from $\mathbb{Q}^a_p$ by topological completion (see, for example, \cite{gouvea}, \cite{robert} for the definition and basic properties of $\mathbb{Q}_p$ and  $\mathbb{C}_p$ ).

 Finally, we use the following standard notation: Given $A,B\subseteq X$ and $x\in X$, we define the best approximation error $E(x,A)=\inf_{a\in A}\|x-a\|$ and the deviation of $B$ from $A$ is given by $E(B,A)=\sup_{b\in B}E(b,A)$. Furthermore, $B(X)=\{x\in X:\|x\|\leq 1\}$ denotes the unit ball of $X$ and, for $r>0$, $B(x,r)=\{y\in X:\|x-y\|\leq r\}$.

\section{An ultrametric lethargy result} \label{section_s}

Let us state the main result of this section:

\begin{theorem}\label{main} Let $(X,\{A_n\})$ be an ultrametric approximation scheme. The following claims are equivalent:
\begin{itemize}
\item[$(a)$] For every sequence of real numbers $\{\varepsilon_n\}\searrow 0$, there exists $x\in X$ such that $E(x,A_n)\neq \mathbf{O}(\varepsilon_n)$.
\item[$(b)$] \[\inf_{n\in\mathbb{N}}E(B(X),A_n)=c>0.\]
\end{itemize}
\end{theorem}

\begin{remark} Let $(X,\|\cdot\|)$ be an ultrametric Banach space. If $X$ is re-normed with an equivalent ultrametric norm $\|\cdot\|^*$, and we denote $B(X)^*=\{x\in X: \|x\|^*\leq 1\}$, $E(x,A_n)^*=\inf_{a\in A_n}\|x-a\|^*$, and  $E(B(X),A_n)^*=\sup_{x\in B(X)^*}E(x,A_n)^*$, we have that
\[\inf_{n\in\mathbb{N}}E(B(X),A_n)>0 \text{ if and only if } \inf_{n\in\mathbb{N}}E(B(X)^*,A_n)^*>0.
\]
This observation is specially useful when $p_{\mathbb{K}}<1$, since every ultrametric
normed space $(X,\|\cdot\|)$ over a discretely valued field $\mathbb{K}$ can be equivalently renormed  with a norm $\|\cdot\|^*$ which satisfies $\|X\|^*=\{\|x\|^*:x\in X\}\subseteq |\mathbb{K}|$ (see, for example, \cite[Theorem 2.1.9]{perez}).
\end{remark}

The proof of Theorem \ref{main} is based on the following three technical results:

\begin{lemma}\label{uno} Let
$h:\mathbb{N}\to\mathbb{N}$ be a map such that $h(n)\geq n$ for all $n$,
 and let $\{\varepsilon _{n}\}\searrow 0$. Then there exists a
sequence
$\{\xi_{n}\}\searrow 0$ such that $\xi_{n}\geq\varepsilon_{n}$ and
$\xi_{n}\leq 2\xi_{h(n)}$ for every $n$.
\end{lemma}

\begin{lemma}\label{cero} Let $\rho\in\mathbb{K}$ with $|\rho|<1$.  Let $x\in X\setminus\{0\}$ and $r\in (0,\infty)$. Then:
 \begin{itemize}
 \item[$(a)$] There exists $\lambda\in \mathbb{K}$ such that
\[|\rho| <\|\lambda x\|\leq 1\]
\item[$(b)$]  There exists $\lambda\in \mathbb{K}$ such that
\[\|\lambda x\|\leq r \text{ and } \frac{1}{|\lambda|}\leq \frac{1}{|\rho|^2}\frac{\|x\|}{r}.\]
\end{itemize}

\end{lemma}


 \begin{lemma}\label{dos} Let $(X,\{A_n\})$ be an ultrametric approximation scheme and assume that  $\{\xi_n\} \searrow 0$  satisfies
 \begin{equation}\label{condicionsalto}
\xi_n\leq C \xi_{K(n+1)-1} \ \ (n\in\mathbb{N})
 \end{equation}
for a certain constant $C>0$. Then the following are equivalent claims:
 \begin{itemize}
 \item[$(a)$] For all $x\in X$ there exists $C(x)>0$ such that \begin{equation}\label{verif_shap} E(x,A_n)\leq C(x)\xi_n, \  \  (n\in\mathbb{N}). \end{equation}
 \item[$(b)$] There exists a constant $D>0$ such that,
 \begin{equation}\label{verif_shap_fuerte} E(x,A_n)\leq D\|x\|\xi_n,\ \ \ (n\in\mathbb{N}, x\in X).\end{equation}
 \end{itemize}
 \end{lemma}
 \noindent \textbf{Proof of Theorem \ref{main}.} Let us prove $(b)\Rightarrow (a)$ or, equivalently, we show that the negation of $(a)$ implies the negation of $(b)$. Assume that $(A_n)$ fails Shapiro's theorem. Then there exists a sequence $\{\varepsilon_n\} \searrow 0$ and a function $C:X\to \mathbb{R}^+$ such that $E(x,A_n)\leq C(x) \varepsilon_n$ for all $n\in \mathbb{N}$ and  all $x\in X$. Then taking the sequence $\{\xi_n\}\searrow 0$ given by Lemma \ref{uno} for the case $h(n)=K(n+1)-1$, and applying Lemma \ref{dos} to this sequence, we have that  (\ref{verif_shap_fuerte}) holds true. It follows that
 \[
 E(B(X),A_n)=\sup_{\|x\|\leq 1}E(x,A_n)\leq D\xi_n \  \ (n\in \mathbb{N}),
 \]
so that $\lim_{n\to\infty}E(B(X),A_n)=0$, which contradicts $(b)$. This ends the proof.

Let us now prove $(a)\Rightarrow (b)$.   Assume that $(b)$ fails. Then  $\lim_{n\to\infty}E(B(X),A_n)=0$. Take $x\in X\setminus\{0\}$ and $\rho\in\mathbb{K}$ such that $0<|\rho|<1$. Then part $(a)$ of Lemma \ref{cero} says us that there exists $\lambda\in \mathbb{K}$ such that $|\rho|<\|\lambda x\|\leq 1$, so that $\frac{1}{|\lambda|}< \frac{1}{|\rho|}\|x\|$. It follows that
\[
E(x,A_n)=\frac{1}{|\lambda|}E(\lambda x,A_n)\leq \frac{1}{|\lambda|}E(B(X),A_n) \leq \frac{1}{|\rho|}\|x\| E(B(X),A_n)
\]
and $(a)$ fails with $\varepsilon_n=E(B(X),A_n)$. This completes the proof of the theorem.  {\hfill $\Box$}



Lemma \ref{uno} is an easy exercise on real sequences:


\noindent \textbf{Proof of Lemma \ref{uno}.} 
Passing from the original function $h$ to
$h^\prime(n) = \max_{1 \leq k \leq n} h(k) + n$, we can assume that
(i) $h(n) > n$ for every $n$, and (ii) the function $h$ is strictly increasing.
Set $m_0 = 0$, and, for $k \geq 1$, $m_k = h(m_{k-1})$. Set $\beta_0 = \varepsilon_1$,
and $\beta_k = \max\{\varepsilon_{m_k}, \beta_{k-1}/2\}$ for $k \geq 1$.
For $n \in \mathbb{N}$, find $k \geq 0$ such that  $n \in [m_k, m_{k+1})$, and set
$\xi_n = \beta_k$.

Then the sequence $(\xi_n)$ has the desired properties.
For $n \in [m_k, m_{k+1})$, $\xi_n = \beta_k \geq \varepsilon_{m_k} \geq \varepsilon_n$.
Furthermore, as $h$ is strictly increasing, $h(n) \in [m_{k+1}, m_{k+2})$, hence
$\xi_{h(n)} = \beta_{k+1} \geq \beta_k/2 = \xi_n/2$.
It remains to show that $\lim \xi_n = 0$, or in other words, that
$\lim \beta_k = 0$. If $\beta_k = \varepsilon_{m_k}$ for infinitely
many values of $k$, then $\lim \beta_k = \lim \varepsilon_{m_k} = 0$.
Otherwise, $\beta_k = \beta_{k-1}/2$ for any $k \geq k_0$. In this case,
too, $\lim \beta_k = 0$.
{\hfill  $\Box$}


\noindent \textbf{Proof of Lemma \ref{cero}.} Let $\rho\in\mathbb{K}$ be such that $0<|\rho|<1$.
Then there exists $a,b\in\mathbb{Z}$ such that $|\rho|^{a+1}< r \leq |\rho|^a$ and $|\rho|^{b+1}< \|x\| \leq |\rho|^b$. Hence
\[
|\rho|=\frac{|\rho|^{b+1}}{|\rho|^{b}}<\frac{\|x\|}{|\rho|^{b}}=\|\frac{x}{\rho^{b}}\|\leq 1,
\]
which proves $(a)$. On the other hand,
\[
\|\rho^{a-b+1}x\|=|\rho|^{a+1}\|\frac{x}{\rho^b}\|\leq r,
\]
and
\[
\frac{1}{|\rho^{a-b+1}|}=|\rho^{b-(a+1)}|=\frac{1}{|\rho|^2}|\rho|^{b+1}\frac{1}{|\rho|^a}\leq \frac{1}{|\rho|^2} \|x\|\frac{1}{r},
\]
which proves part $(b)$ of the lemma. {\hfill $\Box$}

Lemma \ref{dos} is, in our opinion, the difficult part of Theorem \ref{main}.

\noindent \textbf{Proof of Lemma \ref{dos}. } Obviously $(b)\Rightarrow (a)$ (just take $C(x)=D\|x\|$). Let us prove $(a)\Rightarrow (b)$. If $(a)$ holds, then $X=\bigcup_{m\in\mathbb{N}}\Gamma_m$, where
\[
\Gamma_m=\{x\in X: E(x,A_n)\leq m\xi_n, \ \ (n\in\mathbb{N})\}.
\]
The sets $\Gamma_m$ satisfy the following nice properties:
\begin{itemize}
\item[$(i)$] $\Gamma_m=\overline{\Gamma_m}^X$ for all $m\in\mathbb{N}$.
\item[$(ii)$] $\Gamma_m=-\Gamma_m$
\item[$(iii)$]$\Gamma_m+\Gamma_m \subseteq \Gamma_{([C]+1)m}$ for all $m\in\mathbb{N}$, where $C$ is the constant such that $\xi_n\leq C\xi_{K(n+1)-1}$, ($n\in\mathbb{N}$).
\end{itemize}
Property $(i)$ follows from the continuity of the functions $f_n(x)=E(x,A_n)$ and $(ii)$ is a direct consequence of $E(-x,A_n)=E(x,A_n)$. To prove $(iii)$, let us take $x,y\in\Gamma_m$. Then
\[
E(x+y,A_{K(n)})\leq \max\{E(x,A_n),E(y,A_n)\}\leq m\xi_n \ \ (n\in\mathbb{N}).
\]
Thus, given $j\in\mathbb{N}$, $K(n)\leq j\leq K(n+1)-1$, we have that
\[
E(x+y,A_j)\leq E(x+y,A_{K(n)})\leq m\xi_n\leq mC\xi_{K(n+1)-1}\leq m([C]+1)\xi_j.
\]
This proves $\Gamma_m+\Gamma_m\subseteq \Gamma_{([C]+1)m}$.

Taking into account Baire category theorem, there exists $m_0\in\mathbb{N}$ such that $\Gamma_{m_0}$ contains a ball $B(x_0,r_0)$ for a certain $x_0\in X$ and $r_0>0$. Then $B(-x_0,r_0)=-B(x_0,r_0)\subseteq \Gamma_{m_0}$ and, thanks to $(iii)$,
\[
B(0,r_0)=-x_0+B(x_0,r_0)\subseteq -\Gamma_{m_0}+\Gamma_{m_0}=\Gamma_{m_0}+\Gamma_{m_0} \subseteq \Gamma_{([C]+1)m_0}.
\]
Take $x\in X$, $x\neq 0$. If $\rho\in\mathbb{K}$ satisfies $|\rho|<1$, part $(b)$ of Lemma \ref{cero} implies that there exists $\lambda\in\mathbb{K}$ such that $\|\lambda x\|\leq r_0$  and $\frac{1}{|\lambda|}\leq \frac{1}{|\rho|^2}\frac{\|x\|}{r_0}$. Hence $\lambda x\in \Gamma_{([C]+1)m_0}$ and
\[
E(x,A_n)=\frac{1}{|\lambda|}E(\lambda x,A_n)\leq \frac{1}{|\rho|^2}\frac{\|x\|}{r_0} ([C]+1)m_0 \xi_n \ \ (n\in\mathbb{N}, x\in X),
\]
which proves $(b)$.    {\hfill $\Box$}


Once we have proved Theorem \ref{main}, several questions arise. The first one is: Do we have some nice example where the conditions of Theorem \ref{main} are satisfied? Another question is: What can we say about constant $c$ appearing in part $(b)$ of Theorem \ref{main}? We solve both questions with the results below.

\begin{theorem}[Ultrametric Shapiro's Theorem] \label{teosha} If the ultrametric approximation scheme $(X,\{A_n\})$ is such that $A_n$ is a closed linear subspace of $X$ for all $n$, then \[\inf_{n\in\mathbb{N}}E(B(X),A_n)\geq p_{\mathbb{K}}>0.\]
\end{theorem}
\noindent \textbf{Proof. } We just need to prove that, if $Y$ is a proper  closed  subspace of $X$ then there exists $x\in X$ such that $\|x\|\leq 1$ and $E(x,Y)\geq p_{\mathbb{K}}$. Let $y\in X\setminus Y$ then $E(y,Y)=d>0$ since $Y$ is closed. Let $\varepsilon\in (0,1)$ and let $z\in Y$ be such that $\|y-z\|\leq d (1+\frac{\varepsilon}{1-\varepsilon})$.
If $p_{\mathbb{K}}<1$,
Let us take $\rho\in\mathbb{K}$ such that $0<|\rho|<1$ and apply Lemma \ref{cero} with $x=y-z$. Thus, there exists $\lambda\in\mathbb{K}$ such that $|\rho|<\|\lambda x\|\leq 1$. Furthermore, by construction,
\begin{eqnarray*}
\frac{1}{\|x\|} &=& \frac{1}{\|y-z\|} \\
&\geq&  \frac{1}{d}\frac{1}{1+\frac{\varepsilon}{1-\varepsilon}}\\
&=& \frac{1}{d}(1-\varepsilon)
\end{eqnarray*}
Hence $\lambda x\in B(X)$ and
\[
E(\lambda x,Y)=|\lambda|E(x,Y)=|\lambda|E(y,Y)\geq \frac{|\rho|}{\|x\|} d\geq |\rho|(1-\varepsilon).
\]
This proves that
\begin{equation}\label{nueva}
\inf_{n\in\mathbb{N}}E(B(X),A_n)\geq |\rho|,
\end{equation}
since $\varepsilon$ was arbitrarily small. The theorem follows since relation $(\ref{nueva})$ holds for all $\rho\in\mathbb{K}$ with $0<|\rho|<1$.

{\hfill $\Box$}


\begin{theorem}[Dichotomy] \label{dicotomia} Let  $(X,\|\cdot\|)$ be an ultrametric normed  space.
The following are equivalent claims:
\begin{itemize}
\item[$(a)$] \[\inf_{n\in\mathbb{N}}E(B(X),A_n)=c>0.\]
\item[$(b)$]  \[\inf_{n\in\mathbb{N}}E(B(X),A_n)\geq p_{\mathbb{K}}^3>0.\]
\end{itemize}
Furthermore, if    $E(B(X),A_m)< p_{\mathbb{K}}^3$   for a certain $m\in\mathbb{N}$ and the jump function satisfies $K(n)\leq Cn$ for all $n$ and a certain $C>0$, then
$E(B(X),A_n)=\mathbf{O}( \tau^{k^{\frac{1}{\log_2C}}})$ for a certain $0<\tau<1$. In particular, if $C\leq 2$ then the sequence
$\{E(B(X),A_n)\}$ decays exponentially to zero.
\end{theorem}


\noindent \textbf{Proof. } Take $\rho\in\mathbb{K}$ such that $0<|\rho|<1$. We prove that, if $A_n+A_m\subseteq A_{h(n,m)}$, then
\begin{equation}\label{condiciondicotomiacero}
E(B(X),A_{h(n,m)})\leq \frac{1}{|\rho|^2}E(B(X),A_{n})E(B(X),A_{m}).
\end{equation}

Now, the inequality $(\ref{condiciondicotomiacero})$ implies that
\[
\inf_{n\in\mathbb{N}} E(B(X),A_{m})\geq |\rho|^3
\]
since, if  $E(B(X),A_{m_0})<|\rho|^3$ for a certain $m_0\in\mathbb{N}$, then
\[E(B(X),A_{K(m_0)})\leq \left(\frac{E(B(X),A_{m_0})}{|\rho|}\right)^2\leq |\rho|^4
\]
and, if we denote by $K^s(m)=K(K^{s-1}(m))$, $K^0(m)=m$, $K^1(m)=K(m)$, then
\[
E(B(X),A_{K^s(m_0)})\leq |\rho|^{r_s}, \text{ where } r_0=3 \text{ and } r_{s+1}=2(r_s-1).
\]
Obviously, the sequence $(r_s)$ goes to infinity, so that  $E(B(X),A_{K^s(m_0)})$ decreases exponentially. In particular, due to the non-increasing character of the sequence $E(B(X),A_{n})$, we have that $\lim_{n\to\infty}E(B(X),A_n)=0$. This proves the first part of the theorem, since $\rho$ was an arbitrary element of $\{x\in\mathbb{K}:0<|x|<1\}$.

Let us now assume that $K(n)\leq Cn$. We  claim that $E(B(X),A_n)=\mathbf{O}( |\rho|^{\frac{1}{2} k^{\frac{1}{\log_2C}}})$ .  To prove this, we assume with no loss of generality that $K(1)>1$, $E(B(X),A_{1})<|\rho|^3$, and $C>1$ is a natural number. Then $K^s(1)$ converges to infinity and $K^s(1)\leq C^s$ for all $s$. On the other hand, it is clear that $(r_s)$ satisfies $2^s\leq r_s$, $s=0,1,\cdots$, so that
\[
E(C^s)\leq E(K^s(1))=E(B(X),A_{K^s(1)})\leq |\rho|^{r_s}\leq |\rho|^{2^s}.
\]
Thus, taking $s=\log_Ck$, we have that
\begin{eqnarray*}
E(k)&=& E(C^s)\leq E(C^{[s]})\leq |\rho|^{2^{[s]}}\leq |\rho|^{2^{s-1}}\\
&=& \left(|\rho|^{\frac{1}{2}}\right)^{2^{\frac{\log_2k}{\log_2C}}}\\
&=& |\rho|^{\frac{1}{2} k^{\frac{1}{\log_2C}}}
\end{eqnarray*}
Obviously, for $C\leq 2$, $E(k)$ decreases exponentially.

Let us now prove the inequality $(\ref{condiciondicotomiacero})$.
Let $\delta>0$ be fixed and let $x\in X\setminus\{0\}$. Use Lemma \ref{cero} to find $\lambda\in\mathbb{K}$ such that \[|\rho|<\|\lambda x\|\leq 1.\]
Then, by definition of  $E(B(X),A_{n})$, there exists $a_n\in A_n$ such that
\[
\|\lambda x-a_n\|\leq (1+\delta)E(\lambda x, A_n)\leq (1+\delta)E(B(X),A_{n}).
\]
Hence
\[
\| x-\frac{1}{\lambda}a_n\|\leq \frac{1}{|\lambda|}(1+\delta)E(B(X),A_{n})< \frac{\|x\|}{|\rho|}(1+\delta)E(B(X),A_{n}).
\]
It follows that, for all $x\in X$ and all $n\in\mathbb{N}$,
\begin{equation}\label{otra}
E(x,A_n)\leq \frac{\|x\|}{|\rho|}E(B(X),A_{n}),
\end{equation}
since $\delta>0$ was arbitrary. Apply $(\ref{otra})$ to $y=x-\frac{1}{\lambda}a_n$, taking into account that $A_n+A_m\subseteq A_{h(n,m)}$. Then
\[
E(x,A_{h(n,m)})\leq E(y,A_m)\leq \frac{\|y\|}{|\rho|}E(B(X),A_{m})\leq \frac{\frac{\|x\|}{|\rho|}(1+\delta)E(B(X),A_{n})}{|\rho|}E(B(X),A_{m}).
\]
This proves $(\ref{condiciondicotomiacero})$, since $\delta>0$ was arbitrary.

{\hfill $\Box$}

\section{Application to p-adic number theory}
Let us consider $\mathbb{C}_p$ as an ultrametric Banach space over $\mathbb{Q}_p$, and let
$$\mathcal{A}_n=\{\alpha\in \mathbb{C}_p: \exists p(t)\in\mathbb{Q}_p[t], \deg (p)\leq n, \ p(\alpha)=0\}$$
denote the set of $p$-adic algebraic elements of degree $\leq n$. Obviously, $\bigcup\mathcal{A}_n=\mathbb{Q}_p^a$ is the field of $p$-adic algebraic numbers over $\mathbb{Q}_p$, which is a dense subset of $\mathbb{C}_p$. Furthermore, it is well known that $\mathcal{A}_n$ is a closed subset of $\mathbb{C}_p$ for all $n$, and $\lambda \mathcal{A}_n\subseteq \mathcal{A}_n$ for all $\lambda\in\mathbb{Q}_p$ and all $n$ (see \cite[page 130]{robert}). Finally, it is also known that $\mathcal{A}_n+\mathcal{A}_m\subset \mathcal{A}_{nm}$. This proves that $(\mathbb{C}_p,\{\mathcal{A}_n\})$ is an ultrametric approximation scheme with jump function $K(n)=n^2$.
The main result in this section is the following theorem about transcendental $p$-adic numbers:
\begin{theorem}\label{transcendental} Let $\{\varepsilon_n\}\searrow 0$. Then there exists $\alpha\in \mathbb{C}_p\setminus\mathbb{Q}_p^a$ such that
\[
E(\alpha, \mathcal{A}_n)=\inf_{a\in \mathcal{A}_n}|\alpha-a|_p \neq\mathbf{O}(\varepsilon_n).
\]
\end{theorem}
To prove this theorem, we will use Theorem \ref{main}, in conjunction with the following well known result from number theory:
\begin{theorem}[Krasner's Lemma]    Let $\mathbb{K}$ be a non-Archimedean complete valued field of characteristic zero, and let $a$ and $b$ be elements of the algebraic closure of $\mathbb{K}$. Let $a_1=a,a_2,\cdots,a_N$ be the algebraic conjugates of $a$ over $\mathbb{K}$. Suppose that $b$ is closer to $a$ than any other of the conjugates of $a$, i.e., \[
|b-a|<\min_{i=2,3,\cdots,N}|a-a_i|.
\]
Then $K(a)\subseteq K(b)$.
\end{theorem}

This result was initially proved by Ostrowski in 1917 \cite{ostrowski} and rediscovered by Krasner in 1946 \cite{krasner} and it is usually named as Krasner's Lemma \cite[page 178, Theorem 5.7.2]{gouvea}.

\noindent \textbf{Proof of Theorem \ref{transcendental}. } Let $n\in\mathbb{N}$, $n\geq 1$, and let $\sigma_{n}\in\mathbb{C}_p$ be a primitive $p^n$-th root of $1$. Then $\sigma_{n}$ is algebraic of degree $p^n-p^{n-1}$, and its minimal polynomial is given by $g(t)=h(t^{p^{n-1}})$, where $h(t)=t^{p-1}+t^{p-2}+\cdots+t+1$. This is well known, but a direct proof follows by applying Eisenstein's irreducibility criterium to $g(t)$ and taking into account that $t^p-1=(t-1)(t^{p-1}+t^{p-2}+\cdots+t+1)$, so that $$(\sigma_{n}^{p^{n-1}}-1)g(\sigma_{n})=(\sigma_{n}^{p^{n-1}}-1)h(\sigma_{n}^{p^{n-1}})=(\sigma_{n}^{p^{n-1}})^p-1=\sigma_{n}^{p^n}-1=0,$$
which implies that $g(\sigma_{n})=0$. Moreover, it is also well known (\cite[page 107]{robert}) that $$|\sigma_{n}-1|=\frac{1}{p^{\frac{1}{p^n-p^{n-1}}}}.$$ If $\tau_n$ denotes an algebraic conjugate of $\sigma_{n}$, then
\[
|\sigma_{n}-\tau_n|=|\tau_n||\frac{\sigma_{n}}{\tau_n}-1|=\frac{1}{p^{\frac{1}{p^n-p^{n-1}}}}\geq  \frac{1}{p^{\frac{1}{p-1}}}=2c_0>0.
\]
since $\frac{\sigma_{n}}{\tau_n}$ is also a $p^n$-th root of $1$. Moreover, it is easy to check that
\[
\inf_{n\geq 1}\frac{1}{p^{\frac{1}{p^n-p^{n-1}}}}=\frac{1}{p^{\frac{1}{p-1}}}=2c_0>0.
\]
Let us now assume that $\alpha\in \mathcal{A}_n$ and $|\sigma_{n}-\alpha|\leq c_0$. Then Krasner's Lemma implies that $\mathbb{Q}_p(\sigma_{n})\subseteq \mathbb{Q}_p(\alpha)$, which is impossible, since $[\mathbb{Q}_p(\sigma_{n}):\mathbb{Q}_p]=p^n-p^{n-1}>n\geq [\mathbb{Q}_p(\alpha):\mathbb{Q}_p]$. In particular, this implies that $E(\sigma_{n},\mathcal{A}_n)\geq c_0$, so that $\inf_{n\in\mathbb{N}}E(B(\mathbb{C}_p),\mathcal{A}_n)\geq c_0>0$ and we can apply Theorem \ref{main}. This completes the proof.

{\hfill $\Box$}

\section{Acknowledgement}
 The author is very grateful to the anonymous referees of this paper for their many interesting comments. Their advice has greatly improved the readability of this paper.

 \bibliographystyle{amsplain}


\bigskip

\footnotesize{J. M. Almira

Departamento de Matem\'{a}ticas. Universidad de Ja\'{e}n.

E.P.S. Linares,  C/Alfonso X el Sabio, 28

23700 Linares (Ja\'{e}n) Spain

Email: jmalmira@ujaen.es

Phone: (34)+ 953648503

Fax: (34)+ 953648575}

\end{document}